\documentclass[11pt,letterpaper]{article}
\usepackage[T1]{fontenc}
\usepackage[utf8]{inputenc}
\usepackage{lmodern}
\usepackage[margin=1in]{geometry}
\usepackage{amsmath,amssymb,amsthm,mathtools}
\usepackage{microtype}
\usepackage[hidelinks]{hyperref}
\hypersetup{
  pdftitle={A symmetric eigenvalue lower bound for Hadamard products},
  pdfauthor={Shengxuan Luo and Zai Yang}
}

\newtheorem{theorem}{Theorem}[section]
\newtheorem{lemma}[theorem]{Lemma}
\newtheorem{proposition}[theorem]{Proposition}
\newtheorem{corollary}[theorem]{Corollary}
\theoremstyle{definition}
\newtheorem{definition}[theorem]{Definition}
\theoremstyle{remark}

\newcommand{\C}{\mathbb C}
\newcommand{\rank}{\operatorname{rank}}
\newcommand{\tr}{\operatorname{tr}}
\newcommand{\diag}{\operatorname{diag}}
\newcommand{\supp}{\operatorname{supp}}
\newcommand{\FW}{\operatorname{FW}}

\title{A symmetric eigenvalue lower bound for Hadamard products}
\author{%
Shengxuan Luo\thanks{Email:
\href{mailto:2592425519@stu.xjtu.edu.cn}{\texttt{2592425519@stu.xjtu.edu.cn}}.}
\qquad
Zai Yang\thanks{Corresponding author. Email:
\href{mailto:yangzai@xjtu.edu.cn}{\texttt{yangzai@xjtu.edu.cn}}.}\\[0.5em]
{\small School of Mathematics and Statistics, Xi'an Jiaotong University}\\
{\small Xi'an 710049, China}
}
\date{}

\begin{document}
\emergencystretch=2em
\maketitle

\begin{abstract}
For a positive semidefinite matrix $X$, let $\mu_k(X)$ be the
minimum of the least eigenvalues of its principal submatrices of
order $k$. For positive semidefinite $A,B$ of order $n$, we prove
\[
 \mu_{i+j-1}(A\circ B)\geq\mu_i(A)\mu_j(B),
 \qquad i,j\in\{1,2,\dots,n\},\quad i+j-1\leq n,
\]
and characterize equality when the right-hand side is positive.
The resulting symmetric least-eigenvalue bound
\[
 \lambda_{\min}(A\circ B)\geq
 \max_{1\leq i\leq n}\mu_i(A)\mu_{n-i+1}(B)
\]
recovers the classical Schur lower bounds at endpoint choices and can be exact even when both factors are singular. It is incomparable with the projection bound
from a recent companion paper, and combining the two bounds can improve on both. The proof combines a spectral minorant of $B$
whose diagonal entries are at least $\mu_j(B)$ with the projection
inequality underlying the companion bound.
Multifactor and Khatri--Rao consequences are also given.
\end{abstract}

\begingroup

\par
\endgroup

\section{Introduction}\label{sec:introduction}

The Hadamard product of matrices $A=[a_{\ell q}]$ and
$B=[b_{\ell q}]$ of the same size is $A\circ B=[a_{\ell q}b_{\ell q}]$.
Schur's product theorem states that
$A\circ B$ is positive semidefinite whenever $A,B$ are positive semidefinite
\cite{Schur1911,HornJohnson2013}. Spectral inequalities for this
operation have a substantial literature; see, for example,
\cite{Styan1973,AndoHornJohnson1987,WangZhang1995,MondPecaric1998,Zhan1997}.
In particular, a Hadamard product can be positive definite even
when both factors are singular. Rank and Kruskal-rank criteria
for this phenomenon appear in
\cite{YangStoicaTang2019,HornYang2020,HornLuoXuYang2026}.

Throughout this paper, matrices are complex, $X^*$ denotes the conjugate
transpose of a matrix $X$, and $X\succeq0$ means that $X$ is Hermitian positive
semidefinite. 
We write $X_{\ell q}$ for the $(\ell,q)$ entry of $X$.
Write $r_X=\rank X$ for the rank of $X$ and $X[S]$
for a principal submatrix of $X$ with both rows and columns indexed by a set $S$.
For a Hermitian matrix $X$ of order $n$, we order its eigenvalues decreasingly:
\[
	\lambda_1(X)\geq\cdots\geq\lambda_n(X)=\lambda_{\min}(X),
\]
and define
\begin{equation}\label{eq:def-mu}
 \mu_k(X)=\min_{\substack{S\subseteq\{1,\dots,n\}\\|S|=k}}
                  \lambda_{\min}(X[S]),\qquad k\in\{1,\dots,n\}.
\end{equation}
Equivalently,
\begin{equation}\label{eq:sparse-rayleigh}
 \mu_k(X)=
 \min_{\substack{x\in\C^n,\ \|x\|=1\\|\supp x|\leq k}}x^*Xx.
\end{equation}
Thus $\mu_k(X)$ is a uniform lower bound on the Rayleigh quotient
over all coordinate subspaces of dimension at most $k$. In particular,
\begin{equation}\label{eq:mu-properties}
 \mu_1(X)\geq\cdots\geq\mu_n(X)=\lambda_{\min}(X),
 \qquad
 \mu_k(X[S])\geq\mu_k(X)\quad (|S|\geq k).
\end{equation}
Both assertions follow directly from the families over which the
minima are taken.

\begin{theorem}[Main inequality]\label{thm:main}
Let $A,B\succeq0$ be of order $n$ and let $i,j\in\{1,2,\dots,n\}$.
If $i+j-1\leq n$, then
\begin{equation}\label{eq:main}
 \mu_{i+j-1}(A\circ B)\geq\mu_i(A)\mu_j(B).
\end{equation}
Consequently,
\begin{equation}\label{eq:optimized}
 \lambda_{\min}(A\circ B)\geq
 \max_{1\leq i\leq n}\mu_i(A)\mu_{n-i+1}(B).
\end{equation}
\end{theorem}

Since $\mu_1(X)$ is the least diagonal entry of $X$ and
$\mu_n(X)=\lambda_{\min}(X)$, the endpoint choices in
\eqref{eq:optimized} recover the classical Schur lower bounds
\cite{HornJohnson2013,HornLuoXuYang2026}:
\[
 \lambda_{\min}(A\circ B)\geq
 \max\left\{
 \mu_1(A)\lambda_{\min}(B),\
 \lambda_{\min}(A)\mu_1(B)
 \right\}.
\]
When both $A$ and $B$ are singular, these endpoint bounds vanish,
whereas intermediate choices of $i$ in \eqref{eq:optimized}
can yield a strictly positive lower bound.

For ordinary products of $A,B\succeq0$ of order $n$ and
$i,j\in\{1,\dots,n\}$, the classical multiplicative inequalities
are \cite[Theorem~7]{MerikoskiKumar2004}
\begin{subequations}\label{eq:ordinary-product}
\begin{alignat}{2}
 \lambda_{i+j-1}(AB)&\leq\lambda_i(A)\lambda_j(B),
 &\qquad &i+j-1\leq n,\label{eq:ordinary-product-upper}\\
 \lambda_{i+j-n}(AB)&\geq\lambda_i(A)\lambda_j(B),
 &\qquad &i+j\geq n+1.\label{eq:ordinary-product-lower}
\end{alignat}
\end{subequations}
Inequality \eqref{eq:main} may be viewed as a reversal in direction
of the upper inequality \eqref{eq:ordinary-product-upper}: it retains
the output index $i+j-1$ and the condition $i+j-1\leq n$, while
replacing the ordinary product by the Hadamard product and
$\lambda_k$ by $\mu_k$. It is also analogous to the lower
inequality \eqref{eq:ordinary-product-lower}, but with output index
$i+j-1$ rather than $i+j-n$. For the least eigenvalue of $AB$,
the latter inequality requires $i=j=n$, yielding
$\lambda_{\min}(AB)\geq\lambda_{\min}(A)\lambda_{\min}(B)$.
In contrast, \eqref{eq:optimized} uses all complementary orders
$i+j=n+1$ and can remain positive when both factors are singular.

The quantities have different variational meanings:
$\lambda_k(X)$ maximizes the least Rayleigh quotient over
$k$-dimensional subspaces, whereas $\mu_k(X)$ minimizes it over
$k$-dimensional coordinate subspaces. Hence
$\mu_k(X)\leq\lambda_k(X)$.

For $B\ne0$, fix an orthonormal eigenbasis of $B$ ordered by
nonincreasing eigenvalues, and let $P_r$ denote the orthogonal
projection onto the span of its first $r$ vectors. The companion paper
\cite[Theorem~2]{HornLuoXuYang2026} gives the {\em directional projection bound}
\begin{equation}\label{eq:projection-bound}
 \lambda_{\min}(A\circ B)\geq L_{\rm proj}(A\mid B)
 :=\max_{1\leq r\leq r_B}
 \left\{\mu_{n-r+1}(A)\lambda_r(B)
                  \min_\ell(P_r)_{\ell\ell}\right\}.
\end{equation}
We show in this paper that the symmetric bound \eqref{eq:optimized} and the larger of the
two directional projection bounds $L_{\rm proj}(A\mid B)$ and $L_{\rm proj}(B\mid A)$ from \eqref{eq:projection-bound}
are incomparable. In addition, applying the symmetric bound and the directional projection bound to separate positive
semidefinite summands of one factor can nevertheless improve on
both; see Section~\ref{sec:comparison}.

The main technical contribution in the proof of Theorem~\ref{thm:main}
is a uniform diagonal estimate for a trace-balanced upper spectral
truncation of $B$: its diagonal entries are at least $\mu_j(B)$.
This truncation is a spectral minorant of $B$ and a nonnegative
combination of rank-$j$ orthogonal projections. We combine the diagonal estimate with the
projection inequality underlying the companion bound
\eqref{eq:projection-bound}. The same estimate also yields the equality
characterization of \eqref{eq:main} and \eqref{eq:optimized} at interior indices.

Section~\ref{sec:proof} proves the main inequality from the
spectral-minorant lemma and treats equality and sharpness.
Section~\ref{sec:minorant-construction} proves the spectral-minorant lemma.
Section~\ref{sec:comparison} compares and combines the two bounds.
Section~\ref{sec:extensions} gives multifactor and Khatri--Rao
consequences, and Section~\ref{sec:conclusions} concludes this paper.

\section{Proof, equality, and sharpness}\label{sec:proof}

\subsection{The spectral minorant and the main proof}

\begin{definition}[Spectral minorants and upper truncation]
\label{def:spectral-truncation}
Let $B\succeq0$ have distinct eigenvalues
$\theta_1>\cdots>\theta_d\geq0$, and let $E_r$ be the orthogonal
projection onto the eigenspace corresponding to $\theta_r$.
The ordered decomposition $B=\sum_{r=1}^d\theta_rE_r$ is uniquely
determined by $B$. A matrix $Z=\sum_{r=1}^d z_rE_r$ with
$0\leq z_r\leq\theta_r$ for every $r$ is called a
\emph{spectral minorant} of $B$. For $t\geq0$, the
\emph{upper spectral truncation} of $B$ at $t$ is
\[
 T_t(B)=\sum_{r=1}^d\min\{\theta_r,t\}E_r.
\]
\end{definition}
In particular, $T_t(B)$ is a spectral minorant of $B$ and satisfies
\[
 0\preceq T_t(B)\preceq B,\qquad T_t(B)\preceq tI.
\]

\begin{lemma}[Spectral minorant]\label{lem:cap}
Let $B\succeq0$ be of order $n$, and suppose that $j\in\{1,2,\dots,n\}$
and $\mu_j(B)>0$. Set
\begin{equation}\label{eq:tau}
 \tau=\max\left\{t\geq0:\tr\bigl(T_t(B)\bigr)\geq jt\right\},
 \qquad Z=T_\tau(B).
\end{equation}
Then $\tau>0$,
\begin{equation}\label{eq:cap}
 \begin{aligned}
 &0\preceq Z\preceq B,\qquad Z\preceq\tau I,\qquad \tr Z=j\tau,\\
 &Z_{\ell\ell}\geq\mu_j(B)\quad(\ell\in\{1,\dots,n\}).
 \end{aligned}
\end{equation}
In particular, $Z$ is a nonnegative linear combination of rank-$j$
orthogonal projections. If $j<n$, then
\begin{equation}\label{eq:cap-equality}
 Z_{\ell\ell}=\mu_j(B)\quad\Longleftrightarrow\quad
 B_{\ell\ell}=\mu_j(B).
\end{equation}
\end{lemma}

We prove Lemma~\ref{lem:cap} in Section~\ref{sec:minorant-construction}.
The other ingredient is the following consequence of
\cite[Theorem~5]{HornLuoXuYang2026}: if $A\succeq0$ is of order
$n$ and $P$ is a rank-$r$ orthogonal projection, then
\begin{equation}\label{eq:projection}
 A\circ P\succeq\mu_{n-r+1}(A)\diag(P),\qquad r\in\{1,\dots,n\},
\end{equation}
where $\diag(P)$ is the diagonal matrix with the same diagonal as
$P$. Indeed, all principal submatrices of $A-\mu_{n-r+1}(A)I$
of order $n-r+1$ are positive semidefinite, so the cited theorem
applies to its Hadamard product with $P$. A short factor-width proof
of \eqref{eq:projection} is also given below.

\begin{proof}[Proof of Theorem~\ref{thm:main}]
If $\mu_i(A)\mu_j(B)=0$, Schur's product theorem proves
\eqref{eq:main}. Otherwise fix $J\subseteq\{1,\dots,n\}$ of order
$k=i+j-1$. By \eqref{eq:mu-properties},
$\mu_i(A[J])\geq\mu_i(A)$ and $\mu_j(B[J])\geq\mu_j(B)>0$.
Apply Lemma~\ref{lem:cap} to $B[J]$ to obtain the upper
spectral truncation $Z=T_\tau(B[J])$.
Write $Z=\sum_\ell t_\ell P_\ell$, where $t_\ell\geq0$ and
each $P_\ell$ is an orthogonal projection of rank exactly $j$.
Since $k-j+1=i$, apply \eqref{eq:projection} to each $P_\ell$
and use $\mu_i(A[J])\geq\mu_i(A)$ to obtain
\begin{equation}\label{eq:main-matrix}\begin{split}
 (A\circ B)[J]
 & \succeq A[J]\circ Z\\ 
 &=\sum_\ell t_\ell\bigl(A[J]\circ P_\ell\bigr)\\
 &\succeq\sum_\ell t_\ell\mu_i(A)\diag(P_\ell)\\
 &=\mu_i(A)\diag(Z)\\
 &\succeq\mu_i(A)\mu_j(B)I_k.
\end{split}\end{equation}
The first inequality uses Schur's product theorem and
$Z\preceq B[J]$. Taking the minimum over $J$ proves \eqref{eq:main}.
Taking $j=n-i+1$ and maximizing over $i$ gives \eqref{eq:optimized}.
\end{proof}

\subsection{Equality cases}

\begin{proposition}[Equality at interior indices]\label{prop:equality}
Let $A,B\succeq0$ be of order $n$. Suppose that
$i,j\in\{2,\dots,n\}$, $i+j-1\leq n$, and $\mu_i(A)\mu_j(B)>0$.
Equality holds in \eqref{eq:main} if and only if there is a coordinate
$\ell$ such that
\begin{equation}\label{eq:interior-equality}
 A_{\ell\ell}=\mu_i(A),\qquad B_{\ell\ell}=\mu_j(B).
\end{equation}
Each such coordinate is isolated in both factors: all off-diagonal
entries in its row and column vanish.
\end{proposition}
\begin{proof}
Suppose equality holds. Choose a set $J$ of order $k=i+j-1$
and a unit vector $x$ in the coordinates of $J$ attaining the
least eigenvalue $\mu_i(A)\mu_j(B)$ of $(A\circ B)[J]$.
By \eqref{eq:main-matrix}, the diagonal bound for the minorant
of $B[J]$ and the positivity of $\mu_i(A)$ give
\[
 \sum_{\ell\in J}|x_\ell|^2\bigl(Z_{\ell\ell}-\mu_j(B)\bigr)=0.
\]
Each summand is nonnegative, so $Z_{\ell\ell}=\mu_j(B)$ for
every $\ell$ in the support of $x$.
As $Z_{\ell\ell}\geq\mu_j(B[J])\geq\mu_j(B)$, both inequalities
are equalities there. Since $j<k$, \eqref{eq:cap-equality} gives
$B_{\ell\ell}=\mu_j(B)$. Apply the same argument with the factors
interchanged and with the same $J$ and $x$. Since $i<k$, it yields
$A_{\ell\ell}=\mu_i(A)$ on the support of $x$ as well.

Conversely, a coordinate satisfying \eqref{eq:interior-equality}
gives the Rayleigh quotient $\mu_i(A)\mu_j(B)$ on any $k$-set
containing it. Together with \eqref{eq:main} this proves equality.
Finally, as $i\geq2$, every $2\times2$ principal submatrix of
$A-\mu_i(A)I$ is positive semidefinite. A zero diagonal entry
in such a block forces the associated off-diagonal entry to vanish.
The argument for $B$ is identical.
\end{proof}

If either factor has no isolated coordinate, every interior
inequality \eqref{eq:main} with a positive right-hand side is
therefore strict.

We next consider equality in \eqref{eq:main} when an index is $1$
or the right-hand side vanishes.

Suppose first that $i=1$ and $\mu_1(A)\mu_j(B)>0$.
Equality in \eqref{eq:main} holds if and only if there exist
$J\subseteq\{1,\dots,n\}$ with $|J|=j$ and a nonzero vector $x$ indexed
by $J$ such that
\begin{equation}\label{eq:boundary-equality}
 \begin{gathered}
 \bigl[A[J]\circ(B[J]-\mu_j(B)I_j)\bigr]x=0,\\
 x_\ell=0\quad
       (\ell\in J,\ A_{\ell\ell}>\mu_1(A)).
 \end{gathered}
\end{equation}
Indeed, for every such $J$,
\[
 \begin{aligned}
 &(A\circ B)[J]-\mu_1(A)\mu_j(B)I_j\\
 &\quad=A[J]\circ\bigl(B[J]-\mu_j(B)I_j\bigr)\\
 &\qquad+\mu_j(B)
       \bigl(\diag(A[J])-\mu_1(A)I_j\bigr).
 \end{aligned}
\]
Both summands are positive semidefinite. Thus the least eigenvalue
of $(A\circ B)[J]$ equals $\mu_1(A)\mu_j(B)$ precisely when the
two summands have a common nonzero null vector. Minimizing over
$J$ gives the criterion above. The case $j=1$ follows by
interchanging the factors.

Boundary equality need not involve an isolated coordinate.
For $n\geq2$, let $\mathbf1\in\C^n$ be the all-ones vector.
Then $A=\mathbf1\mathbf1^*$ and $B=I+\mathbf1\mathbf1^*$
have no isolated coordinates, yet
$\lambda_{\min}(A\circ B)=\mu_1(A)\lambda_{\min}(B)=1$.

Finally, if $\mu_i(A)\mu_j(B)=0$, equality in
\eqref{eq:main} holds precisely when some principal submatrix
of $A\circ B$ of order $i+j-1$ is singular.
A vanishing right-hand side alone does not imply equality.

Equality in \eqref{eq:optimized} holds precisely when equality
in \eqref{eq:main} holds for a maximizing complementary pair
$i+j=n+1$. If the maximum is positive and attained at an
interior index $i\in\{2,\dots,n-1\}$,
Proposition~\ref{prop:equality} therefore characterizes equality
without any singularity assumption on the factors. If a positive
maximum is attained only at endpoint indices, the boundary
criterion \eqref{eq:boundary-equality} applies.

In particular, if both factors are singular and the right-hand side
of \eqref{eq:optimized} is positive, equality holds precisely when
they have a common minimum-diagonal coordinate and some complementary
orders $i+j=n+1$ satisfy $\mu_i(A)=\mu_1(A)$ and
$\mu_j(B)=\mu_1(B)$. Indeed, both endpoint terms vanish, so
Proposition~\ref{prop:equality} applies to a maximizing pair.
Conversely, such a pair yields $\mu_1(A)\mu_1(B)$, which is
also the product's diagonal entry at the common coordinate.

\subsection{Sharpness of the output index}
\begin{proposition}\label{prop:index}
If $i,j\in\{1,\dots,n\}$ and $i+j-1<n$, there exist positive
semidefinite $A,B$ such that
\[
 \mu_i(A)>0,\qquad \mu_j(B)>0,\qquad
 \mu_{i+j-1}(A\circ B)>0,\qquad \mu_{i+j}(A\circ B)=0.
\]
Thus the index in Theorem~\ref{thm:main} is optimal.
\end{proposition}
\begin{proof}
For $k\geq0$, put $v_k=(1^k,2^k,\dots,n^k)^T$, and take
\[
 A=\sum_{k=0}^{i-1}v_kv_k^*,\qquad B=\sum_{k=0}^{j-1}v_kv_k^*.
\]
The restrictions of $v_0,\dots,v_{r-1}$ to any $r$ coordinates
form a Vandermonde basis. In particular, $\mu_i(A),\mu_j(B)>0$.
Since $v_p\circ v_q=v_{p+q}$,
\begin{equation}\label{eq:vandermonde}
 A\circ B=\sum_{p=0}^{i-1}\sum_{q=0}^{j-1}v_{p+q}v_{p+q}^*.
\end{equation}
Every $v_0,\dots,v_{i+j-2}$ occurs in the sum. By the same
Vandermonde argument, every principal submatrix of order $i+j-1$
is positive definite. As $r_{A\circ B}=i+j-1$, every principal
submatrix of order $i+j$ is singular.
\end{proof}

\section{Proof of the spectral-minorant lemma}
\label{sec:minorant-construction}

This section proves Lemma~\ref{lem:cap}.
The lemma supplies a minorant of $B$ to which the projection
inequality \eqref{eq:projection} can be applied in the proof
of Theorem~\ref{thm:main}. For this purpose, the minorant must
be a nonnegative combination of rank-$j$ orthogonal projections.
We obtain this projection structure by truncating the eigenvalues
of $B$ from above at a threshold chosen through a trace balance. A separate
issue is whether this spectral modification retains the diagonal
lower bound required in Lemma~\ref{lem:cap}. Establishing that
bound is the main task of this section. We first introduce a
sparse test-matrix principle, then apply it to the truncation, and
finally prove the decomposition used in that application.
Throughout the truncation construction, $B\succeq0$, $j\in\{1,\dots,n\}$, and
$\mu_j(B)>0$.

\subsection{The sparse test-matrix principle}

For $s\in\{1,\dots,n\}$, the Hermitian factor-width cone consists of
finite nonnegative linear combinations of the indicated rank-one
matrices:
\begin{equation}\label{eq:fw}
 \FW_s(n)=\operatorname{cone}\{zz^*:z\in\C^n,\ |\supp z|\leq s\};
\end{equation}
see \cite{BomanChenParekhToledo2005}.
The matrix itself may be dense; the support restriction is on
its generating vectors.

By \eqref{eq:sparse-rayleigh}, every vector $z$ supported on at
most $j$ coordinates satisfies
\[
 z^*Bz\geq\mu_j(B)\|z\|^2.
\]
For $M\in\FW_j(n)$, write $M=\sum_\nu z_\nu z_\nu^*$ with
$|\supp z_\nu|\leq j$, absorbing the nonnegative coefficients
into the vectors. Summing the preceding inequalities and using
$\tr(BM)=\sum_\nu z_\nu^*Bz_\nu$ and
$\tr M=\sum_\nu\|z_\nu\|^2$ gives
\begin{equation}\label{eq:local-trace-test}
 \tr(BM)\geq\mu_j(B)\tr M,\qquad M\in\FW_j(n).
\end{equation}
This implication uses only the sparse factors, not the duality
theory developed below.

It therefore suffices to construct, for each coordinate $\ell$,
a test matrix $M_\ell$ satisfying
\begin{equation}\label{eq:test-matrix-goal}
 M_\ell\in\FW_j(n),\qquad
 \tr(BM_\ell)=Z_{\ell\ell},\qquad \tr M_\ell\geq1.
\end{equation}
Indeed, these properties immediately give
\[
 Z_{\ell\ell}=\tr(BM_\ell)
 \geq\mu_j(B)\tr M_\ell\geq\mu_j(B).
\]
We next derive $M_\ell$ from the upper spectral truncation. The technical work in the
last two subsections justifies its factor-width membership.

\subsection{The trace-balanced truncation and its diagonal estimate}

The only decomposition needed in the truncation argument is the following.

\begin{lemma}[Projection-column decomposition]\label{lem:column}
Let $P$ be an orthogonal projection of rank $r$ on $\C^n$.
If $h\geq1$ is an integer with $r+h\leq n$, then for every
$\ell\in\{1,\dots,n\}$,
\begin{equation}\label{eq:column}
 (I-P)e_\ell e_\ell^*(I-P)
 +\frac{P_{\ell\ell}}h(I-P)\in\FW_{r+h}(n),
\end{equation}
where $e_\ell$ denotes the $\ell$th standard coordinate vector in
$\C^n$.
\end{lemma}

The proof of Lemma~\ref{lem:column} is given in
Subsection~\ref{subsec:column-proof}. We first use it to complete
the spectral-minorant construction.

\begin{proof}[Proof of Lemma~\ref{lem:cap}]
\emph{Trace balance.}
Since $\mu_j(B)>0$, we have $r_B\geq j$. The feasible set
in \eqref{eq:tau} is closed, bounded by $\tr B/j$, and contains
$\lambda_j(B)>0$. At its largest element $\tau$, the inequality
must be an equality: if $\tr T_\tau(B)>j\tau$, continuity
would allow a slightly larger feasible threshold. Hence
$\tr Z=j\tau$. In particular, $\tau\geq\lambda_j(B)>0$,
and the order relations in \eqref{eq:cap} follow directly from
Definition~\ref{def:spectral-truncation}.

\emph{Projection representation.}
The matrix $Z/\tau$ has eigenvalues in $[0,1]$ summing to $j$.
Such eigenvalue vectors form the convex hull of the zero-one
vectors with exactly $j$ ones: a vector with any fractional
entries has at least two, which can be perturbed in opposite
directions, preserving their sum, so it is not an extreme point.
Using the same eigenbasis now writes $Z/\tau$ as a convex
combination of rank-$j$ projections. Multiplication by $\tau$
gives the required nonnegative combination for $Z$; see also
\cite{OvertonWomersley1992}.

It remains to establish the diagonal estimate and its equality
condition.
Let $P$ be the projection onto the eigenspaces
of $B$ associated with eigenvalues strictly greater than $\tau$.
Since $\tau\geq\lambda_j(B)$, its rank $r$ is at
most $j-1$. We have
\begin{equation}\label{eq:cap-split}
 Z=\tau P+B(I-P),\qquad
 \tr\bigl(B(I-P)\bigr)=(j-r)\tau.
\end{equation}

\emph{Derivation of the test matrix.}
The two contributions to $Z_{\ell\ell}$ are
$[B(I-P)]_{\ell\ell}$ and $\tau P_{\ell\ell}$.
Since $P$ is a spectral projection of $B$, it commutes with $B$,
and $(I-P)^2=I-P$. Therefore,
\[
 \tr\bigl(B(I-P)e_\ell e_\ell^*(I-P)\bigr)
 =e_\ell^*(I-P)B(I-P)e_\ell
 =[B(I-P)]_{\ell\ell}.
\]
To reproduce the remaining contribution, use the trace balance
in \eqref{eq:cap-split}:
\[
 \frac{P_{\ell\ell}}{j-r}\tr\bigl(B(I-P)\bigr)
 =\tau P_{\ell\ell}.
\]
These identities motivate the definition
\[
 M_\ell=(I-P)e_\ell e_\ell^*(I-P)
       +\frac{P_{\ell\ell}}{j-r}(I-P).
\]
Indeed,
\[
 \begin{aligned}
 \tr(BM_\ell)
 &=[B(I-P)]_{\ell\ell}
 +\frac{P_{\ell\ell}}{j-r}\tr\bigl(B(I-P)\bigr)\\
 &=[B(I-P)]_{\ell\ell}+\tau P_{\ell\ell}
 =Z_{\ell\ell}.
 \end{aligned}
\]

\emph{Factor width.}
Apply Lemma~\ref{lem:column} with $h=j-r\geq1$.
It gives $M_\ell\in\FW_j(n)$ because $r+h=j\leq n$.
Thus the same trace balance that determines the coefficient
$P_{\ell\ell}/(j-r)$ also makes the required factor-width
index equal to $j$.

\emph{Trace and diagonal estimate.}
Since $\tr(I-P)=n-r$, we have
\[
 \tr M_\ell=1-P_{\ell\ell}
   +\frac{P_{\ell\ell}(n-r)}{j-r}
 =1+\frac{n-j}{j-r}P_{\ell\ell}\geq1.
\]
All three properties in \eqref{eq:test-matrix-goal} have now been
verified. Applying \eqref{eq:local-trace-test} yields the sharper
diagonal estimate
\begin{equation}\label{eq:cap-diagonal}
 \begin{aligned}
 Z_{\ell\ell}=\tr(BM_\ell)
 &\geq\mu_j(B)\tr M_\ell\\
 &=\mu_j(B)\left(1+\frac{n-j}{j-r}P_{\ell\ell}\right).
 \end{aligned}
\end{equation}
In particular, $Z_{\ell\ell}\geq\mu_j(B)$, as required in
\eqref{eq:cap}.

If $j<n$ and $Z_{\ell\ell}=\mu_j(B)$, then
\eqref{eq:cap-diagonal} forces $P_{\ell\ell}=0$, and hence
$Pe_\ell=0$. By \eqref{eq:cap-split}, $Z_{\ell\ell}=B_{\ell\ell}$.
The converse in \eqref{eq:cap-equality} follows from
$Z\preceq B$ and the diagonal bound.
\end{proof}

\subsection{Duality and the elementary-symmetric-function lemma}

It remains to prove the projection-column decomposition used above.
We first recall a duality criterion and establish a scalar eigenvalue
inequality; together they will certify the required sparse decompositions.

For the duality argument, we use the equivalent representation
of each member of $\FW_s(n)$ as a sum of positive semidefinite
matrices supported on $S\times S$, with one summand for each
$S\subseteq\{1,\dots,n\}$ of cardinality $s$. The cones increase with $s$:
$\FW_1(n)$ consists of the nonnegative diagonal matrices, and
$\FW_n(n)$ is the positive semidefinite cone.

On the real vector space of Hermitian matrices of order $n$, we use
the Frobenius inner product
\[
 \langle H,M\rangle
 =\operatorname{Re}\tr(H^*M)=\tr(HM).
\]
We denote dual cones by the superscript $\mathrm d$, reserving
$*$ for conjugate transpose.

The following duality characterization is given for real symmetric
matrices in \cite[Proposition~3.4]{GouveiaKovacecSaee2022}.
Its proof applies unchanged to Hermitian matrices.

\begin{lemma}[Factor-width duality]\label{lem:fw-dual}
The cone $\FW_s(n)$ is closed, and its dual is
\begin{equation}\label{eq:fw-dual}
 \FW_s(n)^{\mathrm d}
 =\{H=H^*:H[S]\succeq0
       \text{ for all }S\subseteq\{1,\dots,n\}\text{ with }|S|=s\}.
\end{equation}
Consequently, a Hermitian matrix $M$ belongs to $\FW_s(n)$ if and
only if $\tr(HM)\geq0$
for every $H$ on the right of \eqref{eq:fw-dual}.
\end{lemma}

To prove that a structured matrix belongs to $\FW_s(n)$, we use
the dual criterion rather than construct all sparse factors
explicitly. For a dual test matrix $H$, positivity of its principal
submatrices implies
\begin{equation}\label{eq:principal-minor-sums}
 e_k(\lambda(H))=\sum_{\substack{S\subseteq\{1,\dots,n\}\\|S|=k}}\det H[S]\geq0,\qquad k\in\{1,\dots,s\},
\end{equation}
where $\lambda(H)$ is its eigenvalue vector. For
$x=(x_1,\dots,x_m)\in\mathbb R^m$, the $k$th elementary symmetric
function is
$e_k(x)=\sum_{\substack{S\subseteq\{1,\dots,m\}\\|S|=k}}\prod_{\ell\in S}x_\ell$,
with $e_0(x)=1$.

We need to turn these nonnegative sums into a bound on the
smallest eigenvalue and the sum of selected eigenvalues.
The following lemma provides two operations: deleting a coordinate
preserves one fewer degree of nonnegativity, and nonnegativity
through degree $s$ bounds the negative contribution of the smallest
coordinate. Its proof uses interlacing and Newton's inequalities
for real-rooted polynomials.

\begin{lemma}\label{lem:gamma}
Let $x=(x_1,\dots,x_m)\in\mathbb R^m$ satisfy $e_k(x)\geq0$ for
$k\in\{1,\dots,s\}$, where $s\in\{1,\dots,m\}$.
\begin{enumerate}
\item If $s\geq2$, deleting any coordinate gives a vector whose
elementary symmetric functions through degree $s-1$ are nonnegative.
\item One has
\begin{equation}\label{eq:gamma}
 e_1(x)+(s-1)\min_{\ell\in\{1,\dots,m\}}x_\ell\geq0.
\end{equation}
\end{enumerate}
\end{lemma}
\begin{proof}
Put $p(t)=\prod_\ell(t+x_\ell)$. Differentiating $m-s$
times removes the terms involving $e_{s+1},\dots,e_m$:
\[
 p^{(m-s)}(t)=
 \sum_{k=0}^s\frac{(m-k)!}{(s-k)!}e_k(x)t^{s-k}.
\]
The coefficients are nonnegative and the leading coefficient is
positive, so this polynomial has no positive root. All its roots
are real by Rolle's theorem; hence they are nonpositive.

Deleting $x_\ell$ replaces $p$ by
$q(t)=p(t)/(t+x_\ell)$, whose roots interlace those of $p$.
Differentiation preserves interlacing
\cite[Theorem~1.47]{Fisk2008}, so the roots of $q^{(m-s)}$
lie between those of $p^{(m-s)}$ and are nonpositive.
Since $q^{(m-s)}$ has positive leading coefficient, all its
coefficients are nonnegative. These coefficients are positive
multiples of the elementary symmetric functions of the remaining
coordinates through degree $s-1$, proving part 1.
Here interlacing is understood in the weak sense, allowing
repeated roots.

For part 2, the case $s=1$ is the assumption. If $s=m$, the
nonnegative coefficients of $p$ imply $x_\ell\geq0$ for all $\ell$.
Suppose $2\leq s<m$. Replacing $x$ by
$x+\varepsilon(1,\ldots,1)$ shifts the roots of
$p^{(m-s)}$ strictly to the left, so the elementary symmetric
functions of the perturbed vector are positive through degree $s$.
It suffices to prove the claim for these perturbations and then
let $\varepsilon\downarrow0$. Relabeling the perturbed vector
as $x$, we henceforth assume $e_k(x)>0$ for $k\in\{1,\dots,s\}$.

If the smallest coordinate is negative, write it as $-u$, $u>0$,
and delete it to obtain $\eta$. The recurrence
$e_k(\eta)=e_k(x)+u e_{k-1}(\eta)$ gives $e_k(\eta)>0$
through degree $s$ and
$u<e_s(\eta)/e_{s-1}(\eta)$. Newton's inequalities
\cite[Theorem~4.8, especially~(4.3.3)]{Fisk2008}, with the
nonstrict form for repeated roots obtained by continuity, give
\[
 \left(\frac{e_k(\eta)}{\binom{m-1}{k}}\right)^2
 \geq
 \frac{e_{k-1}(\eta)}{\binom{m-1}{k-1}}\,
 \frac{e_{k+1}(\eta)}{\binom{m-1}{k+1}},
 \qquad k\in\{1,\dots,s-1\}.
\]
Since the terms through degree $s$ are positive, division shows
that the successive ratios of $e_k(\eta)/\binom{m-1}{k}$ are
nonincreasing. Comparing the ratio at degree $s$ with the first
ratio $e_1(\eta)/(m-1)$, and using
$\binom{m-1}{s}/\binom{m-1}{s-1}=(m-s)/s$, gives
\[
 u<\frac{e_s(\eta)}{e_{s-1}(\eta)}
 \leq\frac{m-s}{s(m-1)}e_1(\eta)
 =\frac{m-s}{s(m-1)}\bigl(e_1(x)+u\bigr).
\]
Rearrangement gives
$(s-1)u\leq(m-s)e_1(x)/m\leq e_1(x)$,
which proves \eqref{eq:gamma}. If the smallest coordinate is
nonnegative, the conclusion is immediate. Passing to the limit
proves part 2.
\end{proof}

\subsection{Proof of the projection-column decomposition}
\label{subsec:column-proof}

We first combine Lemma~\ref{lem:gamma} with factor-width duality
to obtain a decomposition for a projection plus a rank-one term.
This will supply the lower-dimensional decomposition used in the
proof of Lemma~\ref{lem:column}.

\begin{lemma}\label{lem:spike}
Let $R$ be a rank-$r$ orthogonal projection on $\C^n$, and let
\(w\in\mathbb C^n\) be a unit vector. For every $h\in\{1,\dots,r\}$,
\begin{equation}\label{eq:spike}
 R+(h-1)ww^*\in\FW_{n-r+h}(n).
\end{equation}
In particular, every rank-$r$ orthogonal projection, $r\geq1$,
belongs to $\FW_{n-r+1}(n)$.
\end{lemma}
\begin{proof}
Take a Hermitian $H$ whose principal submatrices of order
$n-r+h$ are positive semidefinite, as in \eqref{eq:fw-dual}.
All principal submatrices of smaller orders are also positive
semidefinite, so \eqref{eq:principal-minor-sums} applies through
degree $n-r+h$. Delete the largest $n-r$ eigenvalues.
Each deletion reduces the guaranteed degree by one, by
Lemma~\ref{lem:gamma}(1). The remaining $r$ eigenvalues therefore
have nonnegative elementary symmetric functions through degree $h$.
Part 2 of that lemma gives
\[
 \sum_{k=n-r+1}^n\lambda_k(H)+(h-1)\lambda_n(H)\geq0.
\]

For completeness, the minimum-trace inequality follows directly
from an orthonormal eigenbasis $u_1,\dots,u_n$ of $H$. Since $R$
is a rank-$r$ orthogonal projection,
\[
 \tr(HR)=\sum_{k=1}^n\lambda_k(H)\|Ru_k\|^2,\qquad
 0\leq\|Ru_k\|^2\leq1,\qquad
 \sum_{k=1}^n\|Ru_k\|^2=r.
\]
Among weights in $[0,1]$ with sum $r$, the weighted sum of the
eigenvalues is minimized by assigning unit weight to the $r$
smallest eigenvalues. This and the Rayleigh principle give
\[
 \tr(HR)\geq\sum_{k=n-r+1}^n\lambda_k(H),
 \qquad w^*Hw\geq\lambda_n(H).
\]
Since $h-1\geq0$, adding these inequalities shows that
\[
 \tr\bigl(H[R+(h-1)ww^*]\bigr)\geq0.
\]
This holds for every dual test matrix $H$, so
Lemma~\ref{lem:fw-dual} proves \eqref{eq:spike}.
Taking $h=1$ proves the final assertion.
\end{proof}

We now prove the projection-column decomposition stated in the
truncation argument. The proof applies Lemma~\ref{lem:spike} in dimension
$n-1$ and then lifts each sparse factor by adding at most one
coordinate.

\begin{proof}[Proof of Lemma~\ref{lem:column}]
Coordinate permutations preserve factor width, so we may take
$\ell=1$. Write $p=P_{11}$. If $p=0$, the matrix in
\eqref{eq:column} is $e_1e_1^*$; if $p=1$, it is $(I-P)/h$.
Both have the required factor width, the latter because
Lemma~\ref{lem:spike} gives $I-P\in\FW_{r+1}(n)$.

Assume $0<p<1$. Write the first column of $P$ as
\[
 Pe_1=\begin{pmatrix}p\\ \sqrt{p(1-p)}\,w\end{pmatrix},
 \qquad \|w\|=1.
\]
The vectors in the range of $P$ orthogonal to $Pe_1$ have zero
first coordinate. Let $R$ be the projection onto the orthogonal complement
of their trailing coordinates in $\C^{n-1}$. Then $R$ has
rank $n-r$. Orthogonality to $Pe_1$ also makes these trailing
coordinates orthogonal to $w$, and hence $Rw=w$.
The contribution of $Pe_1/\sqrt p$
to the trailing block of $P$ is $(1-p)ww^*$, so the trailing
block of $I-P$ is $R-(1-p)ww^*$.
Since $h\leq n-r$, Lemma~\ref{lem:spike}, applied in dimension
$n-1$ to the rank-$(n-r)$ projection $R$, gives
\[
 R+(h-1)ww^*\in\FW_{r+h-1}(n-1).
\]
Here the factor-width index is
$(n-1)-(n-r)+h=r+h-1$.
The same cone contains $R$, by the case $h=1$ and nesting.

Set
\[
 F=\begin{pmatrix}-\sqrt{(1-p)/p}\,w^*\\ I_{n-1}\end{pmatrix}.
\]
For a vector in $\ker P$, its first coordinate is
$-\sqrt{(1-p)/p}$ times the inner product of $w$ with its trailing
coordinates. The map $F$ inserts precisely this first coordinate.
Using the block description above, direct multiplication gives
\[
 I-P=F\bigl[R-(1-p)ww^*\bigr]F^*,
 \qquad
 (I-P)e_1=-\sqrt{p(1-p)}\,Fw.
\]
Substitute these two identities into the left-hand side of
\eqref{eq:column}. The coefficient of $ww^*$ inside the lifted
matrix becomes $(h-1)(1-p)$, giving
\begin{equation}\label{eq:lift}
 \begin{aligned}
 &(I-P)e_1e_1^*(I-P)+\frac ph(I-P)\\
 &\qquad=\frac ph
 F\bigl[R+(h-1)(1-p)ww^*\bigr]F^*.
 \end{aligned}
\end{equation}
The bracket equals
\[
 pR+(1-p)\bigl[R+(h-1)ww^*\bigr],
\]
so it belongs to $\FW_{r+h-1}(n-1)$ by convexity. If $z$ has
at most $r+h-1$ nonzero entries, the trailing coordinates of
$Fz$ are exactly $z$, and its first coordinate adds at most
one nonzero entry. Thus applying $F$ and $F^*$ to a sparse
decomposition increases its factor width by at most one.
Together with the nonnegative factor $p/h$, this proves
\eqref{eq:column}.
\end{proof}

As an additional consequence, the case $h=1$ in Lemma~\ref{lem:spike}
also recovers the companion inequality \eqref{eq:projection}:
it puts a rank-$r$ projection $P$
in $\FW_{n-r+1}(n)$. On the support of each generator $zz^*$,
the corresponding principal submatrix of $A-\mu_{n-r+1}(A)I$
is positive semidefinite. Schur's product theorem gives
$(A-\mu_{n-r+1}(A)I)\circ(zz^*)\succeq0$.
Taking the corresponding nonnegative combination gives
\eqref{eq:projection}.

\section{Comparison with the rank-optimized projection bound}
\label{sec:comparison}

\subsection{The projection benchmark and positivity}

We set $L_{\rm proj}(A\mid 0)=0$.
For $B\ne0$, the rank-optimized bound \eqref{eq:projection-bound}
follows from $B\succeq\lambda_r(B)P_r$ and \eqref{eq:projection}.
At an eigenvalue tie, a term whose rank lies inside the tied
block is dominated by the term at the last index of that block:
the eigenvalue is unchanged and both the $\mu$ factor and
the diagonal of the projection can only increase. Thus the
maximum is independent of the basis chosen at ties.
The term $r=r_B$ is the full-range estimate. It dominates
the effective-condition-number bound of
\cite[Corollary~3]{HornLuoXuYang2026}, since
$B\preceq\lambda_1(B)P_{r_B}$.

Write $L_{\rm sym}(A,B)$ for the right-hand side of
\eqref{eq:optimized}. The Kruskal rank $k_X$ of $X$ is the
largest $k$ such that every $k$ columns are linearly independent,
with $k_X=0$ if a column is zero. For $X=U^*U$, the identity
$\ker X=\ker U$ shows that the two matrices have the same column
dependencies. A principal submatrix $X[S]$ is the Gram matrix
of the corresponding columns of $U$. Hence
\[
 \mu_k(X)>0\quad\Longleftrightarrow\quad k\leq k_X.
\]
Consequently,
\begin{equation}\label{eq:symmetric-positive}
 L_{\rm sym}(A,B)>0
 \quad\Longleftrightarrow\quad k_A+k_B\geq n+1.
\end{equation}
Indeed, a positive term requires complementary orders
$i\leq k_A$ and $n-i+1\leq k_B$. The inequality itself does
not require either selected order to equal a rank.

In contrast, the directional condition is
\begin{equation}\label{eq:directional-positive}
 L_{\rm proj}(A\mid B)>0
 \quad\Longleftrightarrow\quad
 k_A+r_B\geq n+1\ \text{and}\ B_{\ell\ell}>0\text{ for all }\ell.
\end{equation}
A positive rank-$r$ term implies $k_A+r\geq n+1$ and a positive
diagonal for $P_r$, hence for $B$. Conversely, the conditions
on the right make the full-range term positive: $B$ and its
range projection have zero diagonal entries at the same
coordinates.
Thus \eqref{eq:symmetric-positive} involves the Kruskal ranks of
both factors and implies \eqref{eq:directional-positive} in
both orientations. The converse can fail.

\subsection{Comparison at a rank index}

If $r_B\geq1$, then
\begin{equation}\label{eq:volume}
 \mu_{r_B}(B)\leq\lambda_{r_B}(B)\min_\ell(P_{r_B})_{\ell\ell}.
\end{equation}
For $\mu_{r_B}(B)=0$ this is immediate. Otherwise apply
Lemma~\ref{lem:cap} with $j=r_B$. The largest feasible threshold
in \eqref{eq:tau} is $\lambda_{r_B}(B)$: at that value all $r_B$
positive eigenvalues contribute the threshold, and at any
larger value at least one contributes strictly less. Hence
$Z=\lambda_{r_B}(B)P_{r_B}$, and \eqref{eq:volume} follows from the
minorant's diagonal bound.

Write a term of $L_{\rm sym}(A,B)$ as $\mu_i(A)\mu_j(B)$ with
$i+j=n+1$. Multiplying \eqref{eq:volume} by $\mu_{n-r_B+1}(A)$ shows that
the term in $L_{\rm sym}(A,B)$ with $j=r_B$ is dominated by the
full-range projection estimate. If this term attains
the symmetric maximum, then $L_{\rm sym}(A,B)\leq L_{\rm proj}(A\mid B)$.
The analogous comparison holds upon interchanging the factors.
Neither comparison applies when $i<r_A$ and $j<r_B$.

The common origin of the two bounds does not imply domination.
The symmetric bound uses a nonnegative combination of rank-$j$
projections, whereas the directional bound uses a single leading
spectral projection of one factor. The projections in the
combination need not be leading spectral projections.

\subsection{Two contrasting examples}

First consider Styan's matrices \cite{Styan1973,HornLuoXuYang2026}:
\begin{equation}\label{eq:example-zero}
 A=\begin{pmatrix}2&1&1\\1&1&0\\1&0&1\end{pmatrix},
 \qquad
 B=\begin{pmatrix}2&1&1\\1&1&1\\1&1&1\end{pmatrix}.
\end{equation}
Here $r_A=r_B=2$, but $k_A=2$ and $k_B=1$. Thus
$L_{\rm sym}(A,B)=0$ by \eqref{eq:symmetric-positive}, and
$L_{\rm proj}(B\mid A)=0$ by \eqref{eq:directional-positive}.
In the other orientation,
\[
 \mu_2(A)=\frac{3-\sqrt5}{2},\qquad
 \lambda_2(B)=2-\sqrt2,\qquad
 \min_\ell(P_2)_{\ell\ell}=\frac12.
\]
The rank-$1$ term vanishes because $\mu_3(A)=0$, so
\[
 L_{\rm proj}(A\mid B)
 =\frac{(3-\sqrt5)(2-\sqrt2)}{4}>0,
 \qquad
 \lambda_{\min}(A\circ B)=\frac{5-\sqrt{17}}{2}.
\]

For the opposite ordering, take the two singular factors
\begin{equation}\label{eq:example-positive}
 A=B=\frac14\oplus Q,\qquad
 Q=I_4-\frac14\mathbf1\mathbf1^*,
\end{equation}
where $\mathbf1\in\C^4$ is the all-ones vector. Since $Q$ is a
rank-$3$ projection and its principal submatrices of order $s$
have least eigenvalue $1-s/4$, we have $r_A=r_B=4$ and
\[
 (\mu_1(A),\ldots,\mu_5(A))=
 \left(\frac14,\frac14,\frac14,0,0\right).
\]
Their Kruskal ranks are therefore $3$. The symmetric bound is
attained at $i=j=3$, and is exact by
Proposition~\ref{prop:equality}, since the first coordinate
attains both selected values. For each $r\leq3$, the leading
rank-$r$ projection has a zero first diagonal entry. The
rank-$4$ projection is $1\oplus Q$ and has minimum diagonal
entry $3/4$. Thus
\begin{equation}\label{eq:example-order}
 L_{\rm proj}(A\mid A)=\frac{3}{64}
 <L_{\rm sym}(A,A)=\frac1{16}
 =\lambda_{\min}(A\circ A).
\end{equation}
This strict improvement occurs with neither selected order
equal to a rank. Order $5$ is minimal for the symmetric bound
to exceed both directional bounds with two singular factors:
if $n\leq4$ and $i<r_A$, $j<r_B$, then $i+j\leq2(n-2)<n+1$.
A positive
complementary-index term must therefore select at least one
rank, and the rank-index comparison applies.

The strict improvement persists without an isolated coordinate. Take
\[
 A=B=\left(\frac14\oplus Q\right)+\frac1{100}vv^*,
 \qquad v=(1,1,-1,1,-1)^T.
\]
Both factors are positive semidefinite with all entries nonzero,
common kernel $\operatorname{span}\{(0,1,1,1,1)^T\}$, and range
projection $1\oplus Q$. Direct calculation gives
\[
 \mu_2(A)=\frac{51-\sqrt{626}}{100},\qquad
 \mu_3(A)=\frac{32-\sqrt{374}}{50},\qquad \mu_4(A)=0.
\]
Thus $i=j=3$ is the only positive complementary-index term.
The projection maximum is attained at rank $4$, giving
\[
 L_{\rm sym}(A,A)=\mu_3(A)^2,\qquad
 L_{\rm proj}(A\mid A)=\frac34\mu_2(A)\lambda_4(A),
\]
where $\lambda_4(A)=(13-\sqrt{61})/20$. Consequently,
\[
 \begin{aligned}
 L_{\rm proj}(A\mid A)&\approx0.05056\\
 &<L_{\rm sym}(A,A)\approx0.06412\\
 &<\lambda_{\min}(A\circ A)\approx0.06760.
 \end{aligned}
\]
Here $\lambda_{\min}(A\circ A)=(419-\sqrt{123482})/1000$.
The symmetric bound therefore strictly exceeds both directional
projection bounds even for dense singular factors.

\subsection{Combining the two bounds}

The two bounds can also be used together.

\begin{proposition}[Combined bound]\label{prop:combined}
Let $A,B\succeq0$ be of order $n$, with $B\ne0$.
For $r\in\{1,\dots,r_B\}$, let $P_r$ be a leading rank-$r$ orthogonal
projection of $B$. For every $t\in[0,\lambda_r(B)]$,
\begin{equation}\label{eq:combined-bound}
 \lambda_{\min}(A\circ B)\geq
 L_{\rm sym}(A,B-tP_r)
 +t\,\mu_{n-r+1}(A)\min_\ell(P_r)_{\ell\ell}.
\end{equation}
\end{proposition}
\begin{proof}
Both summands in $B=(B-tP_r)+tP_r$ are positive semidefinite.
Apply \eqref{eq:optimized} to $A\circ(B-tP_r)$ and
\eqref{eq:projection} to $A\circ(tP_r)$. Adding their
least-eigenvalue lower bounds gives \eqref{eq:combined-bound}.
\end{proof}

Taking $t=0$ recovers the symmetric bound. Taking
$t=\lambda_r(B)$ gives at least the rank-$r$ projection
estimate, since the symmetric term is nonnegative. Maximizing
the right-hand side over the admissible choices in both
orientations therefore gives a bound no smaller than either
original estimate.

Strict improvement is possible even when both factors are
singular. Take
\[
 \begin{gathered}
 P=I_4-\tfrac14\mathbf1\mathbf1^*,\qquad
 u=\tfrac1{\sqrt2}(1,-1,0,0)^T,\\
 Q=P-uu^*,\qquad A=B=P+12Q.
 \end{gathered}
\]
Here $r_A=r_B=3$. The matrix $P$ is the range projection of $A$,
and $Q$ is its leading rank-$2$ projection, with $\lambda_2(A)=13$
and $\min_\ell Q_{\ell\ell}=1/4$.
The $2\times2$ principal submatrices give $\mu_2(A)=1$.
Up to simultaneous row and column permutations, the $3\times3$
principal submatrices have spectra
\[
 \left\{1,\frac{13}{4},13\right\},
 \qquad
 \left\{13,\frac{41-\sqrt{1473}}8,
              \frac{41+\sqrt{1473}}8\right\}.
\]
Consequently,
\[
 \mu_3(A)=\frac{41-\sqrt{1473}}8>\frac14.
\]
The rank-$1$ term in \eqref{eq:projection-bound} vanishes,
the rank-$2$ term is $13\mu_3(A)/4>13/16$, and the rank-$3$
term is $3/4$. Thus
\[
 L_{\rm proj}(A\mid A)=\tfrac{13}{4}\mu_3(A),
 \qquad
 L_{\rm sym}(A,A)=\mu_3(A).
\]
In \eqref{eq:combined-bound}, take $r=2$ and $t=12$, so that
$B-tP_r=P$. Since $\mu_2(P)=1/2$ and $\mu_3(P)=1/4$,
\[
 L_{\rm sym}(A,P)
 =\max\left\{\frac14,\frac{\mu_3(A)}2\right\}
 =\frac14.
\]
Thus the combined bound gives
\[
 \lambda_{\min}(A\circ A)
 \geq \frac14+12\,\mu_3(A)\min_\ell Q_{\ell\ell}
 =\frac14+3\mu_3(A),
\]
which exceeds both original bounds.
Decomposing $A\circ A$ into the symmetric and antisymmetric
subspaces associated with the coordinate pairs $\{1,2\}$ and
$\{3,4\}$ gives $\lambda_{\min}(A\circ A)=13/2$.
Numerically, the symmetric, projection, and exhibited combined
estimates are approximately $0.3275$, $1.065$, and $1.233$,
respectively.

\section{Multifactor and Khatri--Rao consequences}\label{sec:extensions}

\begin{corollary}[Multifactor inequality]\label{cor:multifactor}
Let $A_1,\ldots,A_m\succeq0$ be of order $n$, $m\geq2$, and
let $i_\ell\in\{1,\dots,n\}$ for $\ell\in\{1,\dots,m\}$.
For every $k\in\{1,\dots,n\}$ satisfying
$k+m-1\leq\sum_{\ell=1}^m i_\ell$,
\begin{equation}\label{eq:multifactor}
 \mu_k(A_1\circ\cdots\circ A_m)
 \geq\prod_{\ell=1}^m\mu_{i_\ell}(A_\ell).
\end{equation}
\end{corollary}
\begin{proof}
We prove the assertion for $m\geq1$ by induction. For $m=1$,
the condition is $k\leq i_1$, so the assertion follows from
monotonicity. For the induction step, let $m\geq2$ and put
$d=\max\{1,k-i_m+1\}$. Then $d\in\{1,\dots,n\}$,
$d+m-2\leq\sum_{\ell=1}^{m-1}i_\ell$, and
$d+i_m-1=\max\{k,i_m\}\leq n$.
The induction hypothesis gives
\[
 \mu_d(A_1\circ\cdots\circ A_{m-1})
 \geq\prod_{\ell=1}^{m-1}\mu_{i_\ell}(A_\ell).
\]
Apply Theorem~\ref{thm:main} and monotonicity to
$A_1\circ\cdots\circ A_{m-1}$ and $A_m$ with indices $d,i_m$
to obtain \eqref{eq:multifactor}.
\end{proof}

For a matrix $U$ with $n$ columns, let $U_S$ denote the submatrix
of columns indexed by $S$. The quantities already
defined in \eqref{eq:def-mu} describe their conditioning:
\begin{equation}\label{eq:column-conditioning}
 \sqrt{\mu_k(U^*U)}
 =\min_{|S|=k}\ \min_{\substack{x\in\C^k\\\|x\|=1}}\|U_Sx\|,
 \qquad k\in\{1,\dots,n\}.
\end{equation}
This value is zero if $k$ exceeds the number of rows of $U$.

\begin{corollary}[Khatri--Rao inequality]\label{cor:khatri-rao}
Let $U_1,\dots,U_m$ have $n$ columns each, and let $K$ be their
Khatri--Rao product, formed by taking tensor products of
corresponding columns. With the indices of
Corollary~\ref{cor:multifactor},
\begin{equation}\label{eq:khatri-rao}
 \mu_k(K^*K)
 \geq\prod_{\ell=1}^m\mu_{i_\ell}(U_\ell^*U_\ell).
\end{equation}
\end{corollary}
\begin{proof}
Inner products of tensor products are products of the
corresponding inner products, so
\[
 K^*K=(U_1^*U_1)\circ\cdots\circ(U_m^*U_m).
\]
Applying \eqref{eq:multifactor} completes the proof.
\end{proof}

Taking square roots in \eqref{eq:khatri-rao} and using
\eqref{eq:column-conditioning} gives the corresponding uniform
lower bound over all column subsets of the indicated order for
the given matrices. A positive lower bound uniform over a varying
family follows if the product on the right of \eqref{eq:khatri-rao}
is bounded away from zero.
Since $\mu_k(U^*U)>0$ if and only if $k\leq k_U$, setting
$i_\ell=k_{U_\ell}\geq1$ in \eqref{eq:khatri-rao} shows that
for every $k\in\{1,\dots,n\}$ with
$k+m-1\leq\sum_{\ell=1}^m k_{U_\ell}$,
\begin{equation}\label{eq:kruskal}
 k_K\geq k.
\end{equation}
This qualitative Kruskal-rank inequality is classical
\cite{SidiropoulosBro2000,SidiropoulosLiu2001}.
Inequality \eqref{eq:khatri-rao} supplies a quantitative version
in terms of the least eigenvalues of the column Gram matrices.

An exhaustive evaluation of $\mu_k(X)$ involves $\binom{n}{k}$
principal submatrices, so these bounds need not be cheaper to
compute than the eigenvalue they estimate. Their main use is
structural, or quantitative when uniform principal-submatrix
estimates are available.

\section{Conclusions}\label{sec:conclusions}

We established a symmetric eigenvalue lower bound for
positive semidefinite Hadamard products with a sharp
principal-submatrix index. The spectral minorant also explains
when the bound is exact: at interior indices, positive equality
holds precisely when the two factors share an isolated coordinate
at which their diagonal entries attain the selected values.
The principal-submatrix formulation gives multifactor inequalities
and quantitative Khatri--Rao conditioning bounds.

The symmetric bound and the larger of the two directional projection bounds are incomparable. The latter can be
positive when the symmetric bound vanishes, whereas the
symmetric bound can be strictly larger. Applying the two bounds to separate
positive semidefinite summands of one factor can yield a
strictly larger lower bound than either estimate alone.

\section*{Acknowledgment}
The authors would like to thank Prof.~Roger A.~Horn for his helpful  comments which have improved the quality of the paper.
 
\section*{Declaration of competing interest}
The authors declare that they have no known competing financial
interests or personal relationships that could have appeared to
influence the work reported in this paper.

\section*{Data availability}
No external datasets were used. All illustrative matrices and
exact quantities are specified in the article.

\section*{Declaration of generative AI and AI-assisted technologies}
During the preparation of this work, the authors used generative AI in developing
and checking mathematical arguments, and drafting and revising the text and
LaTeX source. The authors are responsible for the verification of
all mathematical results and references and for the final content
of the manuscript.

\end{document}